\documentclass[conference, 10pt]{IEEEtran}
\IEEEoverridecommandlockouts
\usepackage{todonotes}
\setuptodonotes{fancyline, color=blue!30}
\usepackage[T1]{fontenc} % optional
\usepackage[cmintegrals]{newtxmath}
\usepackage{bm}
\usepackage[numbers, sort&compress]{natbib}
\usepackage{amsmath,amsfonts} %amssymb
\usepackage{tabularx, booktabs}
\newcolumntype{Y}{>{\centering\arraybackslash}X} %tabuarlx, centered
\newcolumntype{Z}{>{\raggedleft\arraybackslash}X} %tabularx, right alignment
\usepackage{algorithmic}
\usepackage{graphicx} %, subcaption
\usepackage[caption=false]{subfig}
\usepackage[export]{adjustbox}
\usepackage{textcomp}
\usepackage{pgfplots}
\usepackage{pgfplotstable} % also loads pgfplots
\usepgfplotslibrary{patchplots, colormaps, groupplots}
\pgfplotsset{compat=1.18}
\usetikzlibrary{arrows, chains, positioning,fit,calc, shapes, shapes.geometric, matrix, decorations.pathmorphing, automata, shapes.gates.logic.US,shapes.gates.logic.IEC,calc, spy, external}
\pgfplotsset{/pgfplots/colormap={winter}{rgb255=(0,0,255) rgb255=(0,255,128)}}

\usepackage{tikz}
\usetikzlibrary{arrows,arrows.meta ,chains, positioning,fit,calc, shapes, shapes.geometric, matrix, external, spy}
\usepackage{ifthen, xstring, calc, pgfopts}
\usepackage[printonlyused, nohyperlinks]{acronym}
\usepgfplotslibrary{patchplots, colormaps, groupplots, dateplot}
\pgfplotsset{compat=1.18}
\usepackage{tikz-uml}
\usepackage{todonotes}
\usepackage{enumitem}
\usepackage{fontawesome}

\tikzset{
 decision/.style  = {diamond, black, draw, node distance=1.6cm},
  process/.style    = {fill=black!5,draw, thick, rounded corners=0.2cm, minimum height = 3em,
  minimum width = 3em, text width=4.2cm, align=center, inner xsep=0.2em, inner ysep=0.2em, node distance=0.5cm},
  object/.style    = {fill=black!5,draw, thick, rectangle, minimum height = 3em,
    minimum width = 3em, align=center, inner sep=1em,node distance=4em},
  pin/.style    = {fill=black!5,draw, thick, rectangle, minimum height = 0.6em,
    minimum width = 0.6em, node distance=-1pt, inner sep=0,font=\relsize{-3.5}},
  start/.style      = {fill=black,draw,circle,node distance=0.5cm}, 
  end/.style = {draw, very thick, circle, fill=white, node distance=0.5cm},
  endcenter/.style = {draw, circle, fill=black, minimum size=1.5em, inner sep=0pt, node distance=4em},
  group/.style      = {color=black,thin,rounded corners=0.8em, rectangle}, 
  groupCaption/.style      = {above=0.2cm,right=0.2cm,fill=white}, 
  input/.style    = {coordinate,node distance=2em}, 
  output/.style   = {coordinate,node distance=2em}, 
  between/.style args={#1 and #2}{ % http://tex.stackexchange.com/a/138828/15602
    at = ($(#1)!0.5!(#2)$)
  }  
}

\usepackage{xcolor}
\definecolor{hsurot}	{cmyk}{.00 1 .59 .26}
\definecolor{hsugrau}	{cmyk}{.38 .37 .39 .15}
\definecolor{hsugelb}	{cmyk}{0 .16 .80 0}
\definecolor{hsublau}	{cmyk}{1 .40 0 .82}
\definecolor{hsuturkis}	{cmyk}{1 .14 .60 .49}
\definecolor{hsugruen}	{cmyk}{.16 .16 .91 .28}
\definecolor{hsubraun}	{cmyk}{.00 .57 1 .17}
\definecolor{hsuorange}	{cmyk}{.01 .87 .77 .13}

\usepackage[ruled,vlined]{algorithm2e}%

\SetCommentSty{mycommfont}
\SetAlCapNameFnt{\small}
\SetAlCapFnt{\small}

\usepackage{mdframed}
\usepackage{amsmath}
\usepackage{caption}
\usepackage{float} % For defining new float types

\newfloat{optimizationmodel}{htbp}{lop}
\floatname{optimizationmodel}{Optimization Model}

\newmdenv[
  backgroundcolor=hsugrau!5, % Set the background color to light grey
  linewidth=1pt, % Set the width of the frame line
  innerleftmargin=5pt, % Margin adjustments
  innerrightmargin=5pt,
  innertopmargin=5pt,
  innerbottommargin=5pt
]{OptimizationModelBox}

\usepackage[printonlyused, nohyperlinks]{acronym}
\acrodef{milp}[MILP]{mixed-integer linear programming}
\acrodef{soc}[SOC]{state of charge}
\acrodef{PID}{Proportional-Integral-Differential}
\acrodef{pla}[PLA]{piecewise linear approximation}
\acrodef{MPC}{model predictive control}
\acrodef{RTO}{Real-Time Optimization}
\acrodef{RE}{renewable energy}
\acrodef{SWO}{Site-wide Optimization}
\acrodef{SoS}{Set of Setpoints}
\usepackage[english]{babel} 	% 'french', 'german', 'spanish', 'danish', etc.
\usepackage{hyperref}
\hypersetup{colorlinks=false, hidelinks}%, unicode=true, linkcolor=[rgb]{0.10,0.05,0.67}, citecolor=[rgb]{0.10,0.05,0.67}, filecolor=[rgb]{0.10,0.05,0.67}, urlcolor=[rgb]{0.10,0.05,0.67}}
\addto\extrasenglish{}

\addto\extrasenglish{}
\addto\extrasenglish{}
\addto\extrasenglish{}
\addto\extrasenglish{}
\addto\extrasenglish{}
\addto\extrasenglish{}
\addto\extrasenglish{}
\begin{document}
\bstctlcite{IEEEexample:BSTcontrol}
% \title{From Static to Dynamic: Adapting Static Optimization Models for Real-Time Flexible Energy Resource Control\\
\title{Combination of Site-Wide and Real-Time Optimization for the Control of Systems of Electrolyzers\\
\thanks{$^\dagger$ These authors contributed equally and are listed alphabetically. \\
The authors gratefully acknowledge the funding of the project eModule (Support code: 03HY116) by the German Federal Ministry of Education and Research, based on a resolution of the German Bundestag. This work is also funded by dtec.bw – Digitalization and Technology Research Center of the Bundeswehr out of the project OptiFlex. dtec.bw is funded by the European Union – NextGenerationEU.}
}
\author{\IEEEauthorblockN{Vincent Henkel$^\dagger$, Lukas P. Wagner$^\dagger$, and Felix Gehlhoff}
\IEEEauthorblockA{\textit{Institute of Automation Technology}\\
Helmut Schmidt University Hamburg, Germany \\
Email:\{vincent.henkel, lukas.wagner, felix.gehlhoff\}@hsu-hh.de
}
\and
\IEEEauthorblockN{Alexander Fay}
\IEEEauthorblockA{\textit{Chair of Automation} \\
Ruhr University, Bochum, Germany\\
Email: alexander.fay@rub.de}} 
\maketitle

\begin{abstract}
The rapid expansion of renewable energy sources has introduced significant volatility and unpredictability in the energy supply chain, necessitating advanced control strategies to ensure grid stability and reliability. Green hydrogen production via electrolysis offers a viable solution for converting and storing this volatile renewable energy. However, the inherent fluctuations of renewable energy sources present challenges for consistent utilization and integration of green hydrogen. This work proposes a two-stage optimization approach, combining site-wide optimization and real-time optimization for managing systems of electrolyzers. By adapting an existing static optimization model, dual use is achieved in both site-wide optimization and real-time optimization. The hierarchical optimization structure, characterized by distinct temporal resolutions, enables effective responses to both dynamic changes and long-term trends. The side-wide optimization layer generates long-term plans based on forecast data, while the real-time optimization layer refines these plans in real-time, accommodating immediate fluctuations and ensuring efficient operation. The results from the case study on a system of electrolyzers demonstrate the method's effectiveness in aligning electrolyzer operation with actual availability of renewable energy. This approach offers a robust framework for optimizing the operation of electrolyzers but also other types of flexible energy resources, contributing to sustainable and economically viable energy management.
\end{abstract}

\newcommand*\emptycirc[1][1ex]{\tikz\draw[thick] (0,0) circle (#1);} 
\newcommand*\halfcirc[1][1ex]{%
  \begin{tikzpicture}
  \draw[fill] (0,0)-- (90:#1) arc (90:270:#1) -- cycle ;
  \draw[thick] (0,0) circle (#1);
  \end{tikzpicture}}
\newcommand*\fullcirc[1][1ex]{%
  \begin{tikzpicture}
  \draw[thick, fill] (0,0) circle (#1);
  \end{tikzpicture}}
\newcommand*\notapplicable[1][1.75ex]{\tikz\draw[fill, hsugrau] (0,0) rectangle (#1,#1);}
\section{Introduction} \label{sec:intro}
The rapid expansion of renewable energy (\acsu{RE}) sources significantly increases the volatility and unpredictability in the energy supply chain \cite{eure}, necessitating advanced control strategies to ensure grid stability and reliability \cite{SCF+22}. This emphasizes the crucial role of energy flexibility \cite{Eur18-en}, which is the ability of a resource to modulate its power generation or consumption \cite{UlAn12}. 

Green hydrogen, which can be produced by electrolysis, offers a viable solution for storing volatile \ac{RE} \cite{LKL+23}. Hydrogen serves as a high-density energy carrier for production facilities and provides solutions for energy transportation and storage \cite{ViCa21}. However, leveraging the full potential of green hydrogen production presents challenges due to the inherent fluctuations of renewable energy sources \cite{OUN+22}. These fluctuations introduce uncertainties in the availability and predictability of \ac{RE} supply, complicating the integration and consistent utilization of green hydrogen  \cite{OUN+22}.

To fully capitalize on the benefits of green hydrogen and integrate it effectively into the energy system, advanced optimization strategies are crucial \cite{FPB+21}. The unpredictable nature of \ac{RE} sources necessitates robust control mechanisms to handle these uncertainties \cite{VHS19}. 
However, existing optimization approaches for control often disregard these uncertainties \cite{LGE+23, RWF23}, and instead focus on economic optimization where coarse resolutions are sufficient \cite{WKR+23}. Yet this static optimization is inherently vulnerable to uncertainties \cite{EnHa12}, underscoring the need for finer resolutions to address both renewable integration and uncertainty management \cite{WRK+23}. Achieving such fine resolutions computationally, however, presents significant challenges regarding the timely generation of schedules \cite{WKR+23}.

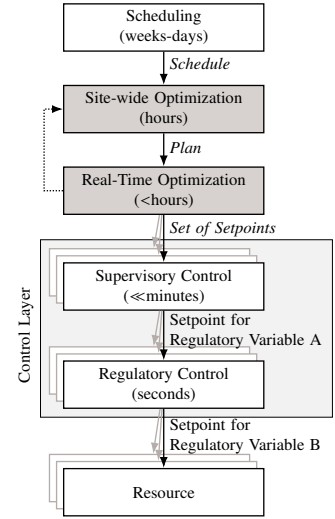
\begin{figure}[ht]
\centering
\resizebox{!}{8cm}{\begin{tikzpicture}[thick, node distance=1.2cm] %font=\sffamily,
    % Styles for processes
    \tikzset{
        process/.style={draw, rectangle, text width = 4cm, minimum height=1cm, minimum width=2.5cm, align=center},
        line/.style={-Latex}
    }

    % Nodes
    \node [process] (scheduling) {Scheduling\\(weeks-days)};
    \node [process, fill=hsugrau!30, below= 0.7cm of scheduling] (siteopt) {Site-wide Optimization\\(hours)};
    \node [process, fill=hsugrau!30, below= 0.7cm of siteopt] (localopt) {Real-Time Optimization\\($<$hours)};
    \node [process,draw=hsugrau!60, fill = white, below= 0.7cm of localopt, xshift = -0.3cm] (supervisory3) {};
    \node [process,draw=hsugrau!60, fill = white, below= 0.85cm of localopt, xshift = -0.15cm] (supervisory2) {};
    \node [process, fill = white, below= 1cm of localopt] (supervisory) {Supervisory Control\\($<\!\!<$minutes)};
    \node [process,draw=hsugrau!60, fill = white, below= 0.75cm of supervisory, xshift = -0.3cm] (regulatory3) {};
    \node [process, draw=hsugrau!60,fill = white, below= 0.9cm of supervisory, xshift = -0.15cm] (regulatory2) {};
    \node [process, fill = white, below= 1.05cm of supervisory] (regulatory) {Regulatory Control\\(seconds)};
    \node[process, draw= none, right = -2.5 cm of supervisory](spacer){};

    % Dotted box
    \begin{scope}[on background layer]
        \node [draw, fill = black!5, fit=(supervisory) (regulatory3) (supervisory3) (regulatory) (spacer), inner sep=5pt] (box) {};
    \end{scope}

    %Text: Distributed Control System 
    \node [align=center,  left = 0.3cm of box, rotate = 90, xshift = 1cm] (pcs) {Control Layer};
    
   \node [process,draw=hsugrau!60, fill = white, below= 0.95cm of regulatory, xshift = -0.3cm] (process2) {};
   \node [process, draw=hsugrau!60, fill = white, below= 1.1cm of regulatory, xshift = -0.15cm] (process3) {};
   \node [process, fill=white, below = 1.25cm of regulatory] (process) {Resource};

    % Lines
    \draw [line] (scheduling) -- (siteopt) node [midway, right, yshift= 0.1cm] {\textcolor{black}{\textit{Schedule}}};;
    \draw [line] (siteopt) -- (localopt) node [midway, right, yshift= 0.05cm] {\textcolor{black}{\textit{Plan}}};
    \draw [line] (localopt) -- (supervisory) node [midway, right, yshift= 0.2cm] {\textcolor{black}{\textit{Set of Setpoints}}};
    \draw [line, draw=hsugrau!60] (localopt) -- (supervisory2);
    \draw [line, draw=hsugrau!60] (localopt) -- (supervisory3);
    \draw [line] (supervisory) -- (regulatory) node [midway, right, yshift= 0.1cm, text width=3.6cm] {Setpoint for\\Regulatory Variable A};
    \draw [line,draw=hsugrau!60] (supervisory) -- (regulatory2);
    \draw [line, draw=hsugrau!60] (supervisory) -- (regulatory3);
    \draw [line] (regulatory) -- (process) node [midway, right, yshift= 0.05cm, text width=3.5cm] {Setpoint for\\Regulatory Variable B};
    \draw [line, draw=hsugrau!60] (regulatory) -- (process2);
    \draw [line, draw=hsugrau!60] (regulatory) -- (process3);
    \draw[densely dotted, line]  (localopt.west) -| node[xshift=-0.7cm, yshift=0.2cm]   {} ++(-.4cm,0cm) |- (siteopt.west) node[near start] {}  node[pos=1] {};

    % Labels
 %   \node [align=center, right= 0.6cm of localopt] (rto) {\ac{RTO}};
 %  \node [align=center, right= 0.6cm of supervisory] (mpc) {\ac{MPC}};
 %   \node [align=center, right= 0.6cm of regulatory] (pid) {\ac{PID}};
\end{tikzpicture}}
\caption{Control Hierarchy (adopted from \cite{Sko23})}
\label{fig:HierarchyOptimizationLayer}
\end{figure}

Therefore, a multi-layered optimization approach, known for its scalability and adaptability, illustrated in \autoref{fig:HierarchyOptimizationLayer}, can be employed. Each layer within this hierarchical structure manages specific decision variables and monitors distinct parameters and variables from other layers. This hierarchical organization facilitates complexity abstraction and enables the provision of services to higher layers. For operational optimization purposes, these layers can be leveraged to achieve various optimization objectives. \cite{Sko23}

A key feature of the hierarchical optimization structure is the distinct temporal resolution of each layer, enabling effective response to both dynamic changes and long-term trends. The \textit{Scheduling} layer, which operates with the longest temporal resolution (weeks to days) \cite{Sko23}, generates a \textit{schedule}. Meanwhile, the \textit{Site-wide Optimization} (\acsu{SWO}) layer, which encompasses a system of flexible energy resources, e.g., multiple electrolyzers, defines a \textit{plan}. Long-term factors, such as market prices, can be integrated by these layers. A \textit{plan} is typically generated for a timespan of hours to days.

The \textit{Real-Time Optimization} (\acsu{RTO}) layer operates at an intermediate resolution, typically hours to minutes\footnote{Please note, that the timescale of \ac{RTO} differs from the shorter timescales seen in other real-time contexts, such as in the communication of controllers.}, aligning operational execution with strategic directives \cite{Sko23}. 
It defines a structured \textit{Set of Setpoints} (\acsu{SoS}), whereby the current setpoint is transferred to the resource via the control layer. The \ac{RTO} layer replans resource operation within one \ac{SWO} time step, taking into account new data, such as short-term forecasts.
Additionally, \ac{RTO}, based on the concept of rolling planning or receding horizon, offers a promising approach to mitigate uncertainties. While traditional strategic optimization methods may disregard this technique, \ac{RTO} enhances adaptability and resilience under uncertainty, making it a powerful tool for economic optimization and improved operational planning \cite{KrSk22}. 

Forecasts, such as those for \ac{RE} generation, are incorporated into these optimization processes. The uncertainty and therefore the accuracy of these forecasts correlates with the forecast horizon \cite{ALW+16}. Consequently, the \textit{schedule} has the highest uncertainty, with decreasing uncertainty down to the \ac{RTO} layer.

The bottom control layers have the shortest temporal resolution, responding to immediate changes within minutes to seconds \cite{Sko23}. In these layers, the overall system is separated into resources, as each resource has its own controller. %, for example.
This layered approach in temporal resolution allows each layer to efficiently address specific tasks while ensuring a coordinated strategy that blends long-term planning with short-term reactions to changed contexts \cite{Sko23}.

The prevailing focus on static optimization \cite{WRK+23} often results in the isolated consideration of the depicted layers. Consequently, many publications disregard the interactions between \ac{RTO} and \ac{SWO} (grey layers in \autoref{fig:HierarchyOptimizationLayer}). %, which inherently complement each other in functionality \cite{Sko23}
Despite its potential for significant economic and energy-systemic benefits, the adoption of \ac{RTO} in real-world applications has fallen short of expectations, leaving its potential largely untapped \cite{KrSk22, EnHa12}.
One possible reason for this can be attributed to the high costs associated with the development and implementation of \ac{RTO} solutions \cite{BaCr08}.

To exploit the untapped potential and ensure more efficient and economic resource operation, this work investigates the transformation of an existing static optimization model (e.g., from \ac{SWO}) into a dynamic \ac{RTO} model. As a result, a two-step optimization approach emerges, utilizing the same optimization model for both long-term and short-term optimization, thus addressing the isolated consideration of the individual layers (see \autoref{fig:HierarchyOptimizationLayer}). The adaptation of the automation pyramid by \citet{Sko23} (as shown in \autoref{fig:HierarchyOptimizationLayer}) is applied in this work. Consequently, the approach is capable of responding to uncertainties and short-term deviations, tackling the issue that different optimization models across layers may lead to inconsistencies \cite{KrSk22}. 
The objectives on the \ac{SWO} and \ac{RTO} layer are somehow similar as both layers focus on, i.e., maximizing benefit, such as minimizing cost \cite{KrSk22}. Therefore the main difference between the two layers is their temporal resolution \cite{KrSk22}. Recently, \citet{RWF23}, have demonstrated that static optimization models can accurately represent the real behavior of a resource, making the models suitable for dual usage in both long-term planning and short-term optimization.

This work is based on previous work by the authors which includes a reusable modular optimization model structure \cite{WRF23b} and a methodology for its parameterization \cite{WaFa24}. 

In summary, the contributions of this work are the following:
\begin{itemize}
    \item A method for adapting existing static optimization models for their dual use in \ac{SWO} and \ac{RTO}
    \item An algorithm for continually solving the two-stage optimization problem incorporating updated forecast information
    \item An evaluation and validation of the proposed method and algorithms through a case study
    \item An assessment of the dynamic \ac{RTO} model on power system performance and efficiency, highlighting improvements over static optimization in terms of cost, stability, flexibility, and energy usage
\end{itemize}

This work is structured as follows: \autoref{sec:relworks} provides an analysis of related work and describes the research gap. \autoref{sec:methodology} describes the method for continually solving the two-stage optimization problem including the approach for adapting existing static optimization models. \autoref{sec:casestudy} evaluates the method via a case study. The method and its evaluation are discussed in \autoref{sec:discussion}. \autoref{sec:conclusion} concludes this work.
\section{Related Work} \label{sec:relworks}
This section reviews related work in optimizing flexible energy resources, focusing on their control, with emphasis on handling uncertainties.

\citet{Sule20} introduced a model for the bidding in energy markets for the integration of renewable energy sources and the operation of flexible energy resources. This approach aims to maximize profits and mitigate risks associated with the uncertainties of renewable energy production. The work by \citet{Sule20} primarily addresses the upper layers of the control hierarchy shown in \autoref{fig:HierarchyOptimizationLayer} -- \textit{Scheduling} and \ac{SWO} -- by proposing a joint bidding strategy for various distributed energy resources based on a predictive model that considers the uncertainties of renewable energy production. However, the model does not address \ac{RTO}, which is crucial for immediate response to short term fluctuations.

\citet{PHM14} describe an optimization model for the operation of a system with energy generation from wind as well as the operation of flexible energy resources under uncertainty. \citet{PHM14} focus on the economic, emission, reliability, and efficiency aspects of the system operation in both certain and uncertain environments of wind, electricity demand, and real-time pricing markets. The authors utilize a stochastic optimization approach, incorporating Monte-Carlo simulations to generate scenario trees based on predicted real-time pricing, wind, and electricity demand. This method allows for the modeling of uncertainties and the identification of optimal operational strategies. Although the generated scenarios address uncertainties, the responsiveness in terms of \ac{RTO} is lacking. This optimization is done on a sub-hourly resolution and cannot be considered as \ac{RTO}, but rather \ac{SWO}.

\citet{VHS19} proposed a model and an algorithm for the operational planning of battery and HVAC resources in buildings using \ac{RTO} to minimize costs of electric energy. \citet{VHS19} also present a two-part forecasting model for short-term variations. This method only focuses on the short-term management for the control of resources without higher planning functions like \ac{SWO}.   This approach underscores the necessity of real-time control but lacks the integration with long-term strategic planning to market system flexibility.%\cite{VHS19}

\citet{TKF+19} proposed a framework for the energy flexible operation of industrial air separation units  using \ac{RTO}. Their work highlights the need to account for process dynamics in production optimization due to the variable nature of electricity prices, and addresses this with a dynamic optimization framework specifically designed for air separation units.
Despite the mention of real-time electricity prices in the study, they are considered at an hourly resolution. Consequently, there is no reaction to short-term deviations, as the optimization primarily focuses on the strategic, longer-term \ac{SWO} of air separation units. 

\citet{FPB+21} presented a specialized optimization model for an electrolyzer. Through experimental analysis of an electrolyzer, detailed linearized models were created to capture the electrolyzer's conversion efficiency and thermal dynamics \cite{FPB+21}. This model informs an \ac{RTO} controller that aims to minimize hydrogen production costs by adapting to fluctuating electricity prices and photovoltaic inflow. 
However, the study by \citet{FPB+21} uses deterministic forecasts, which raises questions about the appropriateness of the method in dealing with forecast uncertainties. Furthermore, \citet{FPB+21} describe their approach as "model predictive control", while \citet{Sko23} highlights the limitations of single-layer optimization in dealing with uncertainty, as the uncertainty must be quantified a-priori. 

\citet{ALY21} introduced an optimization approach for a multi-energy system. Initially, the approach focuses on data management via clustering and scenario reduction to mitigate uncertainties. Subsequently, it employs multi-objective optimization to balance investment and operating costs. The integration of the Markowitz portfolio risk theory allows for managing operational uncertainties. As the approach primarily emphasizes long-term planning and cost optimization, it lacks consideration for real-time resource control to effectively respond to fluctuations, a crucial aspect addressed by \ac{RTO}. Moreover, the hourly resolution used for renewable energy sources proves inadequate in capturing their fluctuations, thus rendering the system not capable of reacting to short-term variations. 

Similar to \citet{ALY21}, \citet{ARA+20} presented a method for the optimization of flexible generation within a smart grid. The method employs linear and non-linear optimization models for solving unit commitment and smart grid scheduling problems, respectively.  However, the method of \citet{ARA+20} does not consider \ac{RTO}, which would be required for the adaptability to short-term variations, but is limited to static, \ac{SWO}. 

\citet{YLS+22} focused on operational optimization for alkaline water electrolysis systems using a mixed-integer nonlinear programming approach.  The optimization considers factors such as solar energy availability, electricity prices, and the resources' operational characteristics for scheduling the electrolysis system. While the optimization effectively plans the operation of the electrolyzers based on the availability of solar energy and electricity prices to increase profitability, it is limited to long-term planning in terms of static, \ac{SWO}. This disregards the aspect of considering short-term fluctuations. As a result, the method cannot fully capture the dynamic nature of renewable energy sources, missing opportunities for further optimization and efficiency improvements in a real-time operational context.

\citet{IrKe24} investigated the management of electric vehicle charging under uncertainty using a two stage optimization approach, thereby focusing on uncertainties such as non-elastic demand and electric vehicle usage behavior. %\cite{IrKe24}

\citet{DDL+21} introduced a hierarchical optimization approach for microgrids, focusing on the coordination between operational planning and \ac{RTO}. At the operational planning level, decisions are made based on day-ahead forecasts to optimize energy costs over one or several days. The \ac{RTO} level adjusts operations based on actual conditions and forecast errors within the current market period. 

Whereas \citet{IrKe24} and \citet{DDL+21} each emphasize the benefits of a two-stage optimization approach for managing flexible energy resources under uncertainty. Instead of reusing existing approaches, both works specialize in applying dedicated optimization models for specific problems within their domains, using different models for the different levels of optimization. However, this can lead to conflicts between the planning and \ac{RTO} layers, potentially resulting in inconsistencies \cite{KrSk22}.

The analysis of related work reveals a gap in the integration of \ac{SWO} and \ac{RTO} into a unified framework, indicating a need for generalized methods to advance the field. Existing approaches focus on either \ac{SWO} or \ac{RTO}, with strategies to address short-term deviations rarely discussed \cite{WRK+23}. This underlines the need to go beyond the use of specific models for individual use cases to develop methods that provide a comprehensive and reusable solution for the dynamic adaptation of optimization strategies.
\section{Method for Solving Two-Stage Optimization} \label{sec:methodology}
This section introduces the method for continually solving \ac{SWO} and \ac{RTO} models for the subsequent control of systems of flexible energy resources. The concept for a two-stage optimization approach using existing optimization models is outlined in \autoref{sec:methodology:concept}. Necessary adjustments for transforming an existing optimization model into a model compatible with both \ac{SWO} and \ac{RTO} are explained in \autoref{sec:method:challenges}. The approach for continually solving the two-stage optimization problem is elaborated in \autoref{sec:method:solving}.

\subsection{Concept for Two-Stage Optimization} \label{sec:methodology:concept}
\begin{figure}[h]
    \centering
    \resizebox{.7\linewidth}{!}{\footnotesize \begin{tikzpicture}[spy using outlines={circle, magnification=6, connect spies}]
%Erster Plot
\begin{axis} [ /pgf/number format/.cd,
        1000 sep={,},
        legend style={
            at={(0,1)},
            anchor=north west,
            nodes={inner sep=0.5pt},
            font=\small, % Hier die Schriftgröße anpassen
        },
        height=3cm, 
        width=5cm,
        xmin=0, xmax=15,
        xtick={0,5,10,15}, % Specify the positions where ticks should be placed
        xticklabels={\(\tau_0\),\(\tau_1\),\(\tau_2\),\(\tau_3\),\(\tau_4\)}, % Manually specify labels
        % minor x tick num=1,
        minor y tick num=1,
        ymin=40, ymax=90, 
        xlabel={Time step $\tau$},
        ylabel={Power},
        ylabel style={xshift=2pt},
        ytick=\empty, 
        %ymajorgrids=true,
        xmajorgrids=true,
        grid style=dashed, 
        name = ax1,
        x tick label style={
            /pgf/number format/.cd,
            fixed,
            fixed zerofill,
            precision=0,
            scaled x ticks=false
        },
        y tick label style={
            /pgf/number format/.cd,
            fixed,
            fixed zerofill,
            precision=0,
        },
    ]

%Schedule
\addplot[mark=,color=hsurot, thick] plot coordinates {
(0,60)
(1,60)
(2,60)
(3,60)
(4,60)
(5,60)
(5,80)
(6,80)
(7,80)
(8,80)
(9,80)
(10,80)
(10,50)
(11,50)
(12,50)
(13,50)
(14,50)
(15,50)
}
;
    % define coordinates at bottom left and top left of rectangle malfunction
  \coordinate (c1) at (axis cs:4.5,85);
  \coordinate (c11) at (axis cs:4.5,45);
  \coordinate (c2) at (axis cs:10.5,45);
  % draw a rectangle
    \draw (c1) rectangle (c2);

\end{axis}

% Zweiter Plot
\begin{axis} [ /pgf/number format/.cd,
        1000 sep={,},
        legend style={
            at={(0.5,-1)},legend columns=-1,
            anchor=south,
            nodes={inner sep=0.5pt},
            font=\footnotesize, % Hier die Schriftgröße anpassen
        },
        height=3cm, 
        width=5cm,
        xmin=0, xmax=11,
        minor x tick num=4,
        minor y tick num=1,
        xtick={0,5,10}, % Specify the positions where ticks should be placed
        xticklabels={\(t_{0,10}\),\(t_{1,5}\),\(t_{1,10}\),\(t_{1,15}\)}, % Manually specify labels
        ymin=35, ymax=100, 
        xlabel={RTO time step $t_{\tau,k}$},
        ylabel={Power},
        ylabel style={xshift=-2pt},
        ytick=\empty, 
        at={($(ax1.south)+(-1.7cm,-2.7cm)$)},
        ymajorgrids=true,
        %xmajorgrids=true,
        grid style=dashed, 
        name = ax2,
        x tick label style={
            /pgf/number format/.cd,
            fixed,
            fixed zerofill,
            precision=0,
            scaled x ticks=false
        },
        y tick label style={
            /pgf/number format/.cd,
            fixed,
            fixed zerofill,
            precision=0,
        },
    ]

% Ist-Schedule
\addplot[mark=*,mark size=1pt,color=hsublau!60, thick, const plot] plot coordinates {
(0,60)
(1,84)
(2,86)
(3,74)
(4,77)
(5,76)
(6,84)
(7,82)
(8,85)
(9,78)
(10,50)
(11,50)
};
\addlegendentry{\textit{SoS}}
% Soll-Schedule
\addplot[mark=,color=hsurot, thick, const plot] plot coordinates {
(0,60)
(1,60)
(1,80)
(2,80)
(3,80)
(4,80)
(5,80)
(6,80)
(7,80)
(8,80)
(9,80)
(10,50)
(11,50)
(12,50)
};
\addlegendentry{\textit{Plan}}
\draw[dashed, hsugrau] (1,100) -- (1,0);
\draw[dashed, hsugrau] (10,100) -- (10,0);
\end{axis}

% draw dashed lines from rectangle in first axis to corners of second malfunction
\draw [dashed] (c2) -- (ax2.north east);
\draw [dashed] (c11) -- (ax2.north west);

%Plan --> RTO
\draw[->, line width=0.9pt, >=Latex]  ([yshift=-0.5cm] ax1.west) -| node[xshift=-0.7cm, yshift=-0.8cm, align=center] {Starting\\ Condition} ++(-.5cm,0cm) |- ([yshift=0.5cm] ax2.west) node[near start] {}  node[pos=1] {};

%RTO --> Plan 
\draw[->, line width=0.9pt, >=Latex]  ([yshift=0.5cm] ax2.east) -| node[xshift=0.6cm, yshift=0.8cm, align=center] {\quad Resource\\\quad Response} ++(.5cm,0cm) |- ([yshift=-0.5cm] ax1.east) node[near start] {}  node[pos=1] {};

% Textfeld für "First Stage Optimization"
\node[align=center, font=\itshape] at ($(ax1.north)+(0cm,0.2cm)$) {Site-wide Optimization};

% Textfeld für "Real Time Optimization"
\node[align=center, font=\itshape] at ($(ax2.north)+(0cm,0.2cm)$) {Real-Time Optimization};

\end{tikzpicture}}    
    % \resizebox{.7\linewidth}{!}{\input{figures/Concept}}
    \caption{Schematic Representation of the Concept}% of the Two-Stage Optimization Approach}
    \label{fig:grundidee}
 \end{figure}
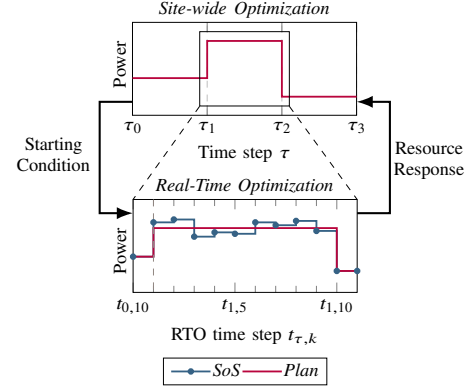

The concept, illustrated in \autoref{fig:grundidee}, utilizes the \ac{SWO} to devise a \textit{plan} for the entire optimization horizon~$\mathcal{T}$, e.g., a full day, segmented into intervals $\Delta \tau$. Individual time steps are denoted as $\tau$. To achieve this, long-term forecast data is incorporated, ensuring the \textit{plan} reflects anticipated future conditions or demands. The results from this planning phase set starting points for \ac{RTO}, i.e., defining starting conditions like resource system states for each time step $t_{\tau}$. This approach is executed cyclically, with new, updated forecasts being taken into account for each optimization.

Following this initial phase, \ac{RTO} refines the optimization at a higher resolution $\Delta t$, such as on a per-minute basis, denoted as time steps $t_{\tau, k}$ over an optimization horizon $T$. The length of the \ac{RTO} horizon $|T|$ is equivalent to one time step of the \ac{SWO} $\Delta \tau$. This finer granularity allows for the accommodation of immediate fluctuations and guarantees that operations can quickly adjust to evolving scenarios. Any changes detected during the \ac{RTO} phase, such as the failure of resources, are then fed back into the \ac{SWO}, ensuring that the behavior of the real system is reflected in the subsequent planning periods $t_{\tau}$ of the \ac{SWO}. This cyclical process creates a dynamic and responsive optimization framework that seamlessly integrates long-term strategic planning with immediate operational adjustments, significantly improving both efficiency and adaptability.

\subsection{Dual Use of Optimization Models} \label{sec:method:challenges}
To facilitate the transformation of static optimization models into dynamic \ac{RTO} models, some preliminary steps must be performed. The basic prerequisite, however, is that a static, feasible optimization model of a system of flexible energy resources is available. The static optimization model requires modification to enable differentiation between its application in \ac{SWO} and \ac{RTO}. Specifically for \ac{RTO} use, it is necessary to fix historic values of variables, including the states of resources or initial setpoints for the \ac{RTO} time steps. This is achieved by equating the respective decision variable to the respective value. These values are derived from the outcomes of the \ac{SWO} or results of previous \ac{RTO} time steps. Furthermore, the energy amount procured for any time step $\tau$ must be set as a target for the corresponding \ac{RTO} horizon starting at $t_{\tau,0}$. This is ensured by means of the constraint shown in \autoref{eq:targetrto}, wherein the energy in one \ac{SWO} time step $\tau$ must be equal to the sum of the energy sourced from the grid in all corresponding \ac{RTO} time steps $t$. 

\begin{align}
    \sum \limits_{t_{{\tau},k} \in T} P_{\text{el, grid},t_{{\tau},k}} \cdot \Delta t = P_{\text{el, grid},\tau} \cdot \Delta \tau \qquad \forall \tau \label{eq:targetrto} % 
\end{align}

This is particularly significant in the context of ancillary services, as the costs for these services are allocated to the parties responsible for deviations from the \textit{plan}.% \cite{}.

Furthermore, the objective function could be modified to ensure the full integration of \ac{RE}, such as maximizing the output of the system.

\subsection{Solving the Two-Stage Optimization Model} \label{sec:method:solving}
\begin{figure}[h]
     \centering
    \resizebox{0.75\linewidth}{!}{\begin{tikzpicture}[auto, thick, >=triangle 45]
   \node [start, minimum width = .65cm ] (start) {};
    % \node [process, below = .5 cm of start] (readModel) {1. Read Optimization Model};
 %   \node [process, below = .5cm of start] (extendModel) {1. Transform Model};
   \node [process,  below = .5cm of start] (initialParameterization) {1. Set Initial Parameters for $\tau_0$};
   \node[rectangle, draw= none, minimum width=1cm, right = 2cm of initialParameterization.east, anchor=east] (optmodel) {\scalebox{3} \faFileTextO};
   \node[rectangle, draw= none, text width=0.8cm, below = 0.3cm of optmodel, text centered, anchor = east] (Textoptmodel) {Optimization \textcolor{white}{ddd}model};
    \node [process, below = .5 cm of initialParameterization] (getForecast) {2. Get Long-Term Forecast for Time Steps $\left[\tau_j, \tau_{|\mathcal{T}|}\right]$};
    \node[rectangle, draw= none, text width=0.8cm, left = 1cm of getForecast.west, text centered, anchor = east] (ForecastIcon) {{\scalebox{3} \faFileTextO}};
    \node[rectangle, draw= none, text width=0.37cm, below = 0.1cm of ForecastIcon, text centered, anchor = east] (TextForecast) {Forecast};
    \node [process, below = .5cm of getForecast] (firstStageOpt) {3. Solve \ac{SWO} Model for $\left[\tau_0, \tau_{|\mathcal{T}|}\right]$};
    \node [process, below = .5cm of firstStageOpt] (saveSchedule) {4. Store \textit{Plan} for $\left[\tau_j, \tau_{|\mathcal{T}|}\right]$};
    \node [process, below = .5cm of saveSchedule] (setStartValues) {5. Initialize RTO with \ac{SWO} Results ($\tau_{j}=t_{\tau_j,0}$)};
    \node [process, below = .5cm of setStartValues] (getShortTermForecast) {6. Get Short-Term Forecast for Time Steps ($\tau_{j}=t_{\tau_j,0}$)};
    \node [process, below = .5cm of getShortTermForecast] (solveModel) {7. Solve \ac{RTO} Model for Horizon $\left[t_{\tau_j,0}, t_{\tau_j,|T|}\right]\subset\tau_j$};
    \node [process, left = 1 cm of solveModel] (setStartValuesRolling) {10. Capture Realized Resource Values for $t_{k}$};
    \node [decision, below = .5cm of  solveModel] (endHorizon) {$t_k=t_{|T|}$};
    \node [right = -0.1cm of endHorizon, text width=1.9cm] {[End of \ac{RTO} horizon?]};
    \node [process, left = 2.35cm of endHorizon] (writeCV) {8. Write Control Variables for $t_k$};
    \node [process,  below = .5cm of writeCV] (getMeasurement) {9. Get Relevant Measurements $t_k$};
    \node [decision, below = .5cm of  endHorizon] (endPlanningHorizon) {$\tau_j=\tau_{|\mathcal{T}|}$};
    \node [left = -0.1cm of endPlanningHorizon, text width=2cm] {[End of Optimization?]};
    \node [end, minimum width = .65cm, below = .5cm of endPlanningHorizon] (end) {};
    \node [draw, minimum width = 4cm, minimum height = 1cm, below = .8cm of getMeasurement] (Process) {Resource};

    \draw[dashed] (initialParameterization) -- (optmodel); 
    \draw [->] (start) -- (initialParameterization);
    \draw [->] (initialParameterization) -- (getForecast);
%    \draw[dashed] (getForecast.west) |- node {} ++(-.5cm,0cm) |- (firstStageOpt.west) node[near start] {}  node[pos=1] {};
%  \draw[dashed] (getForecast.east) |- node {} ++(.5cm,0cm) |- (solveModel.east) node[near start] {}  node[pos=1] {};
    \draw [->] (getForecast) -- (firstStageOpt) ;
    \draw [dashed] (ForecastIcon) ++(0,-1cm) |- (getShortTermForecast);
    \draw [->] (firstStageOpt) -- (saveSchedule);
%    \draw[dashed] (saveSchedule.west) |- node {} ++(-.5cm,0cm) |- (setStartValues.west) node[near start] {}  node[pos=1] {};
    \draw [->] (saveSchedule) -- (setStartValues);
    \draw [->] (setStartValues) -- (getShortTermForecast);
    \draw [->] (getShortTermForecast) -- (solveModel);
    % \draw [->] (solveModel) -- node{$t_k$++} (endHorizon);
        \draw [->] (solveModel) -- (endHorizon);
    \draw [->] (endHorizon) -- node {yes} (endPlanningHorizon);
    \draw [->] (endHorizon.west) -- node {no} (writeCV.east);
    \draw [->] (writeCV) -- (getMeasurement);
    \draw[->]  (getMeasurement.west) |- node {} ++(-.5cm,0cm) |- (setStartValuesRolling.west) node[near start] {}  node[pos=1] {};
    \draw [->] (setStartValuesRolling.east) -- node{$t_k$++} (solveModel.west);
    \draw [->] (endPlanningHorizon) -- node {yes} (end);
    \draw[->]  (endPlanningHorizon.east) -| node[xshift=-0.7cm, yshift=-0.2cm]   {no} ++(2cm,0cm) |- (getForecast.east) node[xshift=1.3cm, yshift=-4.5cm] {$\tau_j$++}  node[pos=1] {};
    \draw [dashed] (getForecast) -- (ForecastIcon);

    \draw[dashed, ->]([xshift=1cm] getMeasurement.south) -- ([xshift=1cm] Process.north);
    \draw[dashed, <-]([xshift=-1cm] getMeasurement.south) -- ([xshift=-1cm] Process.north);

\end{tikzpicture}}
    \caption{Flowchart for the use of an existing optimization model in \ac{RTO}}
    \label{fig:flowchart-rto}
\end{figure}
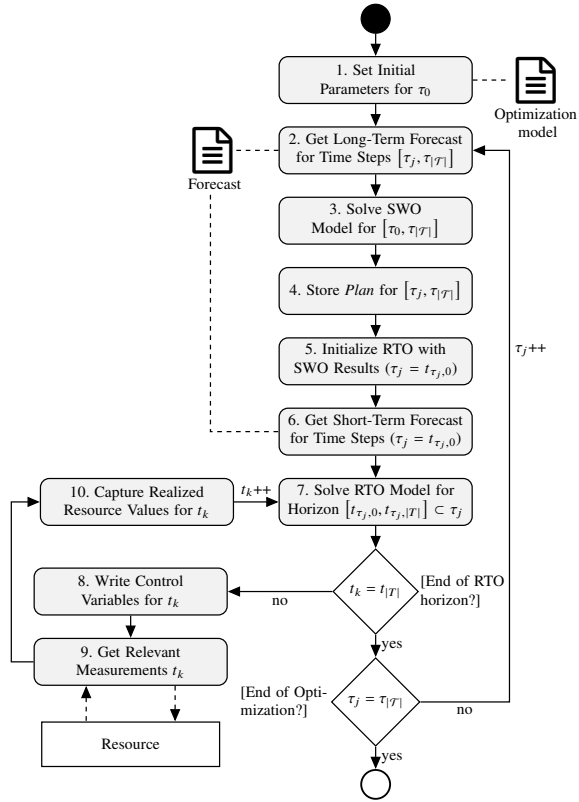

% First, the static optimization model, which can come from the \ac{SWO}, for example, is imported and then modified according to the procedure described in \autoref{sec:method:challenges}, allowing its dual use for both \ac{RTO} and operational planning (Step 1). 
The process for solving the two-stage optimization model is depicted in \autoref{fig:flowchart-rto}. First, initial parameter settings are defined, such as the system state at the start of \ac{SWO} (Step 1).

Subsequently, long-term forecast information, such as electricity data from spot markets or \ac{RE} generation, is imported (Step 2) and utilized to solve the static optimization model for the time steps $\left[\tau_0, \tau_{|\mathcal{T}|}\right]$ (Step~3). The results obtained from this operational planning for the time step $\tau$ are stored (Step~4) and subsequently adopted as initial values for the \ac{RTO} model ($\tau_{j}=t_{\tau_j,0}$, Step~5).

To enable appropriate responses, short-term forecasts with a high resolution, as those from \ac{RE} sources, are incorporated (Step 6). Subsequently, the \ac{RTO} model is solved based on the specified \ac{RTO} horizon $T$, allowing for suitable adjustments to short-term fluctuations (Step~7). Subsequently, the system verifies whether the end of the \ac{RTO} horizon \mbox{$t_k = t_{|T|}$} has been reached. If not, the optimized setpoints are conveyed to the process (Step~8), and relevant measurement values are received (Step~9). These measured values can then be considered during \ac{RTO} (Step~10). The current setpoints then act as the new starting values for the next \ac{RTO} iteration (Step~10). 

Upon reaching the end of the \ac{RTO} horizon with \mbox{$t_k = t_{|T|}$}, the system checks if the end of the planning interval~$\mathcal{T}$ has been reached as well. If so, the process terminates. However, this process is typically executed in a continual manner, meaning that in most cases, the end is never reached. Instead, new starting values for the subsequent planning time step $\tau$ are obtained from the static optimization, and the \ac{RTO} process recommences (Step~5).
Both stages are always solved for the entire length of the respective optimization horizon $\mathcal{T}$/$T$ incorporating newly available information, such as forecasts for future time steps, and historic values. Hence, each model is solved $x$ times, with $x$ corresponding to $\frac{\mathcal{T}}{\Delta \tau}$ or $\frac{T}{\Delta t}$, respectively for \ac{SWO} and \ac{RTO}.

Applying the method outlined in this section generates \textit{plans} for the control of a system of resources. At the availability of new data, such as forecasts, an updated \ac{SoS} or \textit{plan} for the remaining time steps of the respective optimization horizon is generated. %
\section{Evaluation of the Method} \label{sec:validation}  \label{sec:casestudy}
The case study assesses the effectiveness of the proposed method by applying it to a system of electrolyzers, which draws power from both the grid and a wind farm. For this purpose, this section describes the setup of the validation (\autoref{sec:vali:setup}) as well as the results (\autoref{sec:vali:results} \& \autoref{sec:results:method}).

\subsection{Setup of the Evaluation} \label{sec:vali:setup}
The optimization model used for the validation of the method is built using the validated optimization model structure developed by \citet{WRF23b} (see \autoref{optmodel:sitewide}). Therein, a set of constraints for the representation of flexibility features was developed \cite{WRF23b}. This includes constraints for operational boundaries (min./max. values for flows), input-output relationships, as well as system states and related constraints \cite{WRF23b}. The objective function aims at minimizing the cost of electric energy procured from the European intra-day market (\autoref{eq:obj}).

\begin{align}
    \min \text{Cost} = \sum \limits_{t\in \mathcal{T}} P_{\text{grid},t} \cdot \Delta t \cdot c_{\text{el},t}\label{eq:obj}
\end{align}

\begin{optimizationmodel}
\begin{OptimizationModelBox}
\textbf{min} \quad Total cost of energy procured from European intra-day market (\autoref{eq:obj})\\
\textbf{subject to}
\begin{itemize}
\item Power balance to include grid and renewable power and connect electrolyzers
\item Operational boundaries of the system (\autoref{eq:power-boundaries-resource})
\item Operational boundaries of each resource (\autoref{eq:power-boundaries-resource})
  \item Input-Output relationships (\autoref{eq:pla-poweroutput}-\ref{eq:pla-sum1})
  \item Target for energy input (\autoref{eq:target})
  \item System states (\autoref{eq:states-binary}-\ref{eq:powerlimitsbystatemax})
  \item Follower states (\autoref{eq:statesequences})
  \item Holding durations (\autoref{eq:minholdingduration}-\ref{eq:maxholdingduration})
  \item Ramp limits (\autoref{eq:rampmin}-\ref{eq:rampmax})
  % \item Dependencies of the resources
\end{itemize}
\end{OptimizationModelBox}
\vspace{-1em}
\caption{Electrolyzer model used for the case study. Based on model structure developed by \citet{WRF23b}} \label{optmodel:sitewide}
\end{optimizationmodel}

The parameters for each electrolyzer were determined by employing the methodology for parameterizing optimization models developed by \citet{WaFa24}. This creates a data model for the parameterization of the optimization model, based on operational data \cite{WaFa24}. The mathematical modeling and the parameter set are described in the Appendix in detail.

\autoref{optmodel:rto} shows how the existing, feasible \autoref{optmodel:sitewide} is extended for its use in \ac{RTO} as described in \autoref{sec:method:challenges}. Feasibility is ensured through a) a validated model structure \cite{WRF23b} and b) a systematic derivation of suitable parameters \cite{WaFa24}.

\begin{optimizationmodel}
\begin{OptimizationModelBox}
\textbf{max} hydrogen output \textbf{subject to}
\begin{itemize}
\item Constraints of \autoref{optmodel:sitewide}
\item Target to avoid penalty costs (\autoref{eq:targetrto})
\item Fixed Values for past periods (see \autoref{fig:flowchart-rto})
\end{itemize}
\end{OptimizationModelBox}
\vspace{-1em}
\caption{Model used for \ac{RTO}} \label{optmodel:rto}
\end{optimizationmodel}

\begin{figure}[H]
    \centering
    \subfloat[\ac{SWO} time steps $\tau_{0,1,9}$ \label{fig:forecast-firststage}]{\resizebox{!}{3.8cm}{ \begin{tikzpicture}[spy using outlines={circle, magnification=4, connect spies}]
\begin{axis} [ /pgf/number format/.cd,
        % use comma,
        1000 sep={,},
legend style={at={(0.5,-0.45)},anchor=north, legend columns=4},
width=0.49\linewidth, height=4cm, 
  label style={font=\large},
xmin=0, xmax=9, 
xtick={0,1,...,9}, % Specify the positions where ticks should be placed
xticklabels={\(\tau_0\), \(\tau_1\), \(\tau_2\), \(\tau_3\), \(\tau_4\), \(\tau_5\), \(\tau_6\), \(\tau_7\), \(\tau_8\), \(\tau_9\)},
% xmin=1, xmax=10, 
% xtick = {5,5.5,6}, 
% minor x tick num=1,  
% ymin = 0, ymax = 2500, 
 % ytick = {0,40,80},
minor y tick num=1,
xlabel={Time step $\tau$},
ylabel={RE power in kW},
ylabel style={yshift=-3pt},
ymajorgrids=true,
xmajorgrids=true,
grid style = dotted, 
  y tick label style={
        /pgf/number format/.cd,
            fixed,
            fixed zerofill,
            precision=2,
        /tikz/.cd
    },
% tick style={draw=none}
]
\addplot[const plot, mark=, color=blue]  table[y=tau0, x = t]  {figures/data_first-stage-forecast.txt};
\addlegendentry{$\tau_0$}

\addplot[const plot, mark=, color=red]  table[y=tau1, x = t]  {figures/data_first-stage-forecast.txt};
\addlegendentry{$\tau_1$}

% \addplot[const plot, mark=, color=green]  table[y=tau2, x = t]  {figures/data_first-stage-forecast.txt};
% \addlegendentry{$\tau_2$}

% \addplot[const plot, mark=, color=orange]  table[y=tau3, x = t]  {figures/data_first-stage-forecast.txt};
% \addlegendentry{$\tau_3$}

% \addplot[const plot, mark=, color=purple]  table[y=tau4, x = t]  {figures/data_first-stage-forecast.txt};
% \addlegendentry{$\tau_4$}

% \addplot[const plot, mark=, color=brown]  table[y=tau5, x = t]  {figures/data_first-stage-forecast.txt};
% \addlegendentry{$\tau_5$}

% \addplot[const plot, mark=, color=pink]  table[y=tau6, x = t]  {figures/data_first-stage-forecast.txt};
% \addlegendentry{$\tau_6$}

% \addplot[const plot, mark=, color=lime]  table[y=tau7, x = t]  {figures/data_first-stage-forecast.txt};
% \addlegendentry{$\tau_7$}

% \addplot[const plot, mark=, color=cyan]  table[y=tau8, x = t]  {figures/data_first-stage-forecast.txt};
% \addlegendentry{$\tau_8$}

 \addplot[const plot, mark=, thick, densely dotted, color=hsuturkis]  table[y=tau9, x = t]  {figures/data_first-stage-forecast.txt};
\addlegendentry{$\tau_9$}
\node[draw, thick, circle, hsugrau!50!black, minimum width = .43cm] at (axis cs:1.5,0.35) (circ1) {};
% \coordinate (spypoint) at (axis cs:1.5,0.33);
% \coordinate (magnifyglass) at (axis cs:3,0.3);
\end{axis}
% \spy [blue, size=.5cm] on (spypoint) in node[fill=white] at (magnifyglass);
\end{tikzpicture}}}
    \subfloat[\ac{RTO} time steps of $\tau_0$ \label{fig:forecast-rto}]{\resizebox{!}{3.7cm}{ \begin{tikzpicture}
\begin{axis} [ /pgf/number format/.cd,
        % use comma,
        1000 sep={,},
legend style={at={(0.5,-0.5)},anchor=north, legend columns=4},
width=0.49\linewidth, height=4cm, 
  label style={font=\large},
xmin=0, xmax=9, 
xtick={0,1,...,9}, % Specify the positions where ticks should be placed
xticklabels={\(t_0\), \(t_1\), \(t_2\), \(t_3\), \(t_4\), \(t_5\), \(t_6\), \(t_7\), \(t_8\), \(t_9\)},
% xmin=1, xmax=10, 
% xtick = {5,5.5,6}, 
% minor x tick num=1,  
% ymin = 0, ymax = 2500, 
 % ytick = {0,40,80},
minor y tick num=1,
xlabel={RTO time step $t_{\tau_0,k}$},
ylabel={RE power in kW},
ylabel style={yshift=-3pt},
ymajorgrids=true,
xmajorgrids=true,
grid style = dotted, 
  y tick label style={
        /pgf/number format/.cd,
            fixed,
            fixed zerofill,
            precision=2,
        /tikz/.cd
    },
% tick style={draw=none}
]
\addplot[const plot, mark=, color=cyan]  table[y=t0, x = t]  {figures/data_rto1_forecast.txt};
\addlegendentry{$t_{0}$}

\addplot[const plot, mark=, color=purple]  table[y=t1, x = t]  {figures/data_rto1_forecast.txt};
\addlegendentry{$t_{1}$}

% \addplot[const plot, mark=, color=lime!50!black]  table[y=t2, x = t]  {figures/data_rto1_forecast.txt};
% \addlegendentry{$t_{2}$}

% \addplot[const plot, mark=, color=orange]  table[y=t3, x = t]  {figures/data_rto1_forecast.txt};
% \addlegendentry{$t_3$}

% \addplot[const plot, mark=, color=purple]  table[y=t4, x = t]  {figures/data_rto1_forecast.txt};
% \addlegendentry{$t_4$}

% \addplot[const plot, mark=, color=brown]  table[y=t5, x = t]  {figures/data_rto1_forecast.txt};
% \addlegendentry{$t_5$}

% \addplot[const plot, mark=, color=pink]  table[y=t6, x = t]  {figures/data_rto1_forecast.txt};
% \addlegendentry{$t_6$}

% \addplot[const plot, mark=, color=lime]  table[y=t7, x = t]  {figures/data_rto1_forecast.txt};
% \addlegendentry{$t_7$}

% \addplot[const plot, mark=, color=cyan]  table[y=t8, x = t]  {figures/data_rto1_forecast.txt};
% \addlegendentry{$t_8$}

 \addplot[const plot, mark=, thick, densely dotted, color=hsuturkis]  table[y=t9, x = t]  {figures/data_rto1_forecast.txt};
 \addlegendentry{$t_9$}
\end{axis}
\end{tikzpicture}}}
    \caption{Forecast and their uncertainty (data by \citet{ALW+16})}
    \label{fig:forecast}
\end{figure}
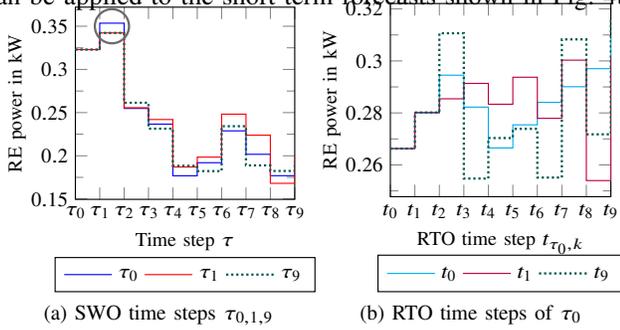

For forecasts of future \ac{RE} generation, this work utilizes data from a real wind farm published by \citet{ALW+16}, available at a resolution of 1 Hz. From this dataset, both $|\mathcal{T}|$ long-term (see \autoref{fig:forecast-firststage}) and $|T|$ corresponding short-term forecasts to each time step~$\tau \in \mathcal{T}$ (see \autoref{fig:forecast-rto}) were generated.

These forecasts were generated in an iterative manner at each time step, retaining data points for past time steps. To address the issue of prediction uncertainties, greater uncertainty for time steps further in the future was introduced. This can be exemplary seen in deviations of values at time step $\tau_1$ of the forecasts generated at $\tau_0$ and $\tau_1$, respectively, in \autoref{fig:forecast-firststage} (grey circle). There, the forecast generated at time step $\tau_0$ (blue) overestimated the renewable generation at $\tau_1$ and necessitated a correction in forecast for $\tau_1$ (red). A similar interpretation can be applied to the short-term forecasts shown in \autoref{fig:forecast-rto}.

The method has been implemented in Java utilizing an optimization model based on IBM ILOG CPLEX\footnote{Implementation: \url{https://github.com/lukas-wagner/TwoStageOpt}}.

\subsection{Results of the Evaluation} \label{sec:vali:results}
\begin{figure}[H]
    \centering
    \resizebox{!}{4cm}{\begin{tikzpicture}
\begin{axis} [ /pgf/number format/.cd,
        % use comma,
        1000 sep={,},
width=0.6\linewidth, height=4cm, 
  label style={font=\large},
xmin=0, xmax=9, 
xtick=data, 
% xtick = {5,5.5,6}, 
% minor x tick num=1,  
% ymin = 0, ymax = 2500, 
ytick = {0,40,60,80},
minor y tick num=1,
xlabel={Time step $\tau$},
ylabel={$c_{\text{el}}$ in EUR/MWh},
 ymajorgrids=true,
 xmajorgrids=true,
grid style = dotted, 
  y tick label style={
        /pgf/number format/.cd,
            fixed,
            fixed zerofill,
            precision=0,
        /tikz/.cd
    },
% tick style={draw=none}
]
\addplot[const plot, mark=, color=hsuturkis, thick]  table[y=price, x = t]  {figures/data_price.txt};
% \addlegendentry{$c_{\text{el}}$}
\end{axis}
\end{tikzpicture}}
    \caption{Electricity price \cite{EPEb}}
    \label{fig:proice}
\end{figure}
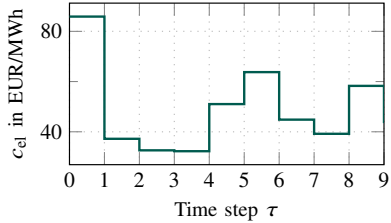

The method outlined in \autoref{sec:methodology} is applied utilizing optimization models~\ref{optmodel:sitewide} and~\ref{optmodel:rto}, with all models solved as depicted in \autoref{fig:flowchart-rto}. Prices sourced from the European intra-day market are employed \cite{EPEb}, as illustrated in \autoref{fig:proice}.

The calculations were performed on Windows 10 with an Intel Core i7-11700 processor and 16 GB RAM with an optimality gap of 10$^{-3}$. 
The method for the two stage optimization approach was applied for $|\mathcal{T}|=10$ \ac{SWO} time steps and $|T|=10$ time steps of \ac{RTO} each. 
The calculation time for each time step $\tau$ including its corresponding \ac{RTO} time steps $t$ was approx. 10~s in total using temporal resolutions of 0.25~h for \ac{SWO} and 0.025~h for \ac{RTO}.

Moreover, the optimized \ac{SoS} generated by \ac{RTO} is transmitted to a simulation model of electrolyzers, as described in \cite{MOCK20241885}, utilizing OPC-UA (Steps 9 \& 10 in \autoref{fig:flowchart-rto}). For this purpose, the method devised by \citet{RWF23} is employed.
This process was carried out successfully and the recorded values are analyzed below.
%\todo[inline]{"The advantages of the two-stage approach are well-explained, but the discussion could benefit from more concrete examples. For instance, the text mentions improved operational planning and economic optimization but does not detail specific scenarios. Providing concrete examples like "In a real-world application, this approach reduced energy procurement costs by 10\%" would illustrate the benefits of the proposed method."}

\subsubsection{Results of Site-wide Optimization} \label{sec:res:sitewide}
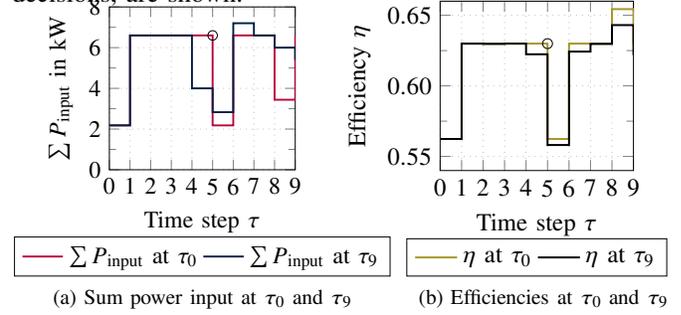
\begin{figure}[H]
\centering
\subfloat[Sum power input at $\tau_{0}$ and $\tau_9$ \label{fig:results:firststage:sumpower}]{\resizebox{!}{3.8cm}{\begin{tikzpicture}
\begin{axis} [ /pgf/number format/.cd,
        % use comma,
        1000 sep={.},
legend style={at={(0.5,-0.45)},anchor=north, legend columns=-1},
width=0.5\linewidth, height=4cm, 
  label style={font=\large},
xmin=0, xmax=9, 
xtick=data, 
% minor x tick num=1,   
ymin = 0, ymax = 8, 
minor y tick num=1,
xlabel={Time step $\tau$},
ylabel={$\sum P_{\text{input}}$ in kW},
% yticklabel pos=left ,axis y line*=left, 
 ymajorgrids=true,
 xmajorgrids=true,
grid style = dotted, 
  y tick label style={
        /pgf/number format/.cd,
            fixed,
            fixed zerofill,
            precision=0,
        /tikz/.cd
    },
% tick style={draw=none}
]
\addplot[const plot, mark=, color=hsurot, thick]  table[y=sum0
, x = t] {figures/data_suminput_efficiency_firststage.txt};
\addlegendentry{$\sum P_{\text{input}}$ at $\tau_0$}
\addplot[const plot, mark=, color=hsublau, thick]  table[y=sum9, x = t] {figures/data_suminput_efficiency_firststage.txt};
\addlegendentry{$\sum P_{\text{input}}$ at $\tau_9$}
\addplot[mark=o, mark size = 2pt, color = black] plot coordinates {(5,6.6)};
% \node[pin= below :{$\tau_4$}] at (axis cs:4,7) {};
\end{axis}
\end{tikzpicture}}}
\subfloat[Efficiencies at $\tau_{0}$ and $\tau_9$ \label{fig:results:firststag:efficiency}]{\resizebox{!}{3.7cm}{\begin{tikzpicture}
\begin{axis} [ /pgf/number format/.cd,
        % use comma,
        1000 sep={,},
legend style={at={(0.5,-0.45)},anchor=north, legend columns=2},
width=0.5\linewidth, height=4cm, 
  label style={font=\large},
xmin=0, xmax=9, 
xtick=data, 
% minor x tick num=1,  
ymin = 0.54, ymax = 0.66, 
minor y tick num=1,
xlabel={Time step $\tau$},
ylabel={Efficiency $\eta$},
% yticklabel pos=right ,axis y line*=right, 
  ymajorgrids=true,
  xmajorgrids=true,
 grid style = dotted, 
  y tick label style={
        /pgf/number format/.cd,
            fixed,
            fixed zerofill,
            precision=2,
        /tikz/.cd
    },
% tick style={draw=none}
]
\addplot[const plot, mark=, color=hsugruen, thick]  table[y=eff0, x = t] {figures/data_suminput_efficiency_firststage.txt};
\addlegendentry{$\eta$ at $\tau_0$}
\addplot[const plot, mark=, color=black, thick]  table[y=eff9, x = t] {figures/data_suminput_efficiency_firststage.txt};
\addlegendentry{$\eta$ at $\tau_9$}
\addplot[mark=o, mark size = 2pt, color = black] plot coordinates {(5,0.629935352)};
\end{axis}
\end{tikzpicture}}}
\caption{Results of \ac{SWO}}
\label{fig:results:firststage}
\end{figure}

\autoref{fig:results:firststage} shows the \textit{plan} generated by the \ac{SWO}, taking into account the forecasts shown in \autoref{fig:forecast-firststage}. For clarity, only the first time step $\tau_0$ and the final, realized \ac{SoS} generated at $\tau_9$, taking into account all past decisions, are shown.

\autoref{fig:results:firststage:sumpower} illustrates that the \textit{plan} generated at $\tau_0$ aligns with the realized \textit{plan} ($\tau_9$) up to and including time step $\tau_4$. This initial \textit{plan} is generated assuming perfect foresight it benefits from the most flexibility potential. % As time progresses, 

Starting from $\tau_5$, a significant share of time steps has already been executed, thereby constraining the ability to respond to updated forecast data during the remaining time steps of \ac{SWO} due to the constraints imposed. Consequently, this results in deviations between the optimized \textit{plan} at $\tau_0$ and realized \textit{plan} from this juncture onwards as the flexibility potential decreases. This limitation reflects in a decrease of 0.5~\% in hydrogen produced between the initial stage at $\tau_0$ and the final step $\tau_9$. The main driver is the uncertainty of the forecasts, with 1.5~\% less realized \ac{RE} than initially forecasted. This also results in a lower overall efficiency of 62.3~\%, compared to the initial \textit{plan} of 62.6~\% at $\tau_0$. The comparison of efficiencies at $\tau_0$ and $\tau_9$ is depicted in \autoref{fig:results:firststag:efficiency}, illustrating the impact of the aforementioned uncertainties on the overall efficiency of the electrolyzer system.

\subsubsection{Results of \textit{Real-Time Optimization}}
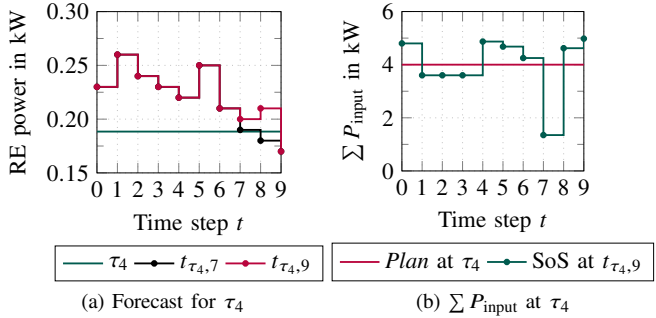
\begin{figure}[H]
\centering
\subfloat[Forecast for $\tau_{4}$ \label{fig:results:rto:ee}]{\resizebox{!}{3.8cm}{\begin{tikzpicture}
\begin{axis}[ /pgf/number format/.cd,
        % use comma,
        1000 sep={,},
legend style={at={(0.5,-0.45)},anchor=north, legend columns=-1},
width=0.49\linewidth, height=4cm, 
  label style={font=\large},
xmin=0, xmax=9, 
xtick=data, 
% minor x tick num=1,    
ymin = 0.15, ymax = 0.3, 
minor y tick num=1,
xlabel={Time step $t$},
ylabel={RE power in kW},
% yticklabel pos=right ,axis y line*=right, 
  ymajorgrids=true,
  xmajorgrids=true,
 grid style = dotted, 
  y tick label style={
        /pgf/number format/.cd,
            fixed,
            fixed zerofill,
            precision=2,
        /tikz/.cd
    },
% tick style={draw=none}
]
\addplot[const plot, mark=, color=hsuturkis, thick]  table[y=eefirst, x = t] {figures/data_rto4_results.txt};
\addlegendentry{$\tau_4$}
\addplot[const plot, mark=*,mark size=1pt, color=black, thick]  table[y=eerto7, x = t] {figures/data_rto4_results.txt};
\addlegendentry{$t_{\tau_4,7}$}
\addplot[const plot, mark=*,mark size=1pt, color=hsurot, thick]  table[y=eerto9, x = t] {figures/data_rto4_results.txt};
\addlegendentry{$t_{\tau_4,9}$}
\end{axis}
\end{tikzpicture}}} %\hspace{-0.4cm} 
\subfloat[$\sum P_{\text{input}}$ at $\tau_{4}$ \label{fig:results:rto:sumpower}]{\resizebox{!}{3.7cm}{\begin{tikzpicture}
\begin{axis} [ /pgf/number format/.cd,
        % use comma,
        1000 sep={,},
legend style={at={(0.5,-0.45)},anchor=north, legend columns=-1},
width=0.49\linewidth, height=4cm, 
  label style={font=\large},
xmin=0, xmax=9, 
xtick=data, 
% minor x tick num=1,  
ymin = 0, ymax = 6, 
minor y tick num=1,
xlabel={Time step $t$},
ylabel={$\sum P_{\text{input}}$ in kW},
% yticklabel pos=right ,axis y line*=right, 
  ymajorgrids=true,
  xmajorgrids=true,
 grid style = dotted, 
  y tick label style={
        /pgf/number format/.cd,
            fixed,
            fixed zerofill,
            precision=0,
        /tikz/.cd
    },
% tick style={draw=none}
]
\addplot[const plot, mark=, color=hsurot, thick]  table[y=suminpfirst, x = t] {figures/data_rto4_results.txt};
\addlegendentry{\textit{Plan} at $\tau_4$}
\addplot[const plot, mark=*,mark size=1pt, color=hsuturkis, thick]  table[y=suminprto9, x = t] {figures/data_rto4_results.txt};
\addlegendentry{\ac{SoS} at $t_{\tau_4,9}$}
\end{axis}
\end{tikzpicture}}}
    \caption{RTO Schedules at $\tau_4$}
    \label{fig:results:rto}
\end{figure}

\autoref{fig:results:rto} shows the results of \ac{RTO} for one exemplary time step $\tau_4$. In \autoref{fig:results:rto:ee}, the prediction for the \ac{RE} output under the \ac{SWO} for time step $\tau_4$, as well as the deviating forecast at $t_{\tau_4,9}$, are depicted. Notably, the output of \ac{RE} sources consistently surpasses the long-term prediction of the \ac{SWO}. Consequently, this leads to an increase in hydrogen production. The realized \ac{SoS} at $t_{\tau_4,9}$ yields 0.66 kWh of hydrogen energy instead of the planned amount of 0.62 kWh at $\tau_4$ capitalizing on the increased availability of \ac{RE} (see \autoref{fig:results:rto:sumpower}). 

Furthermore, as shown in in \autoref{fig:results:rto:sumpower}, a significant drop in power can be observed at $t_7$. This is attributed to the increase in the forecast for \ac{RE} from $t_8$ to $t_9$ compared to previous forecasts (see $t_{\tau_4,7}$). Consequently, the hydrogen production has been shifted to later time steps $t_8$ and $t_9$ resulting in a decrease of power input at $t_7$. 

\subsection{Validation of the Method for Two-Stage Optimization} \label{sec:results:method}
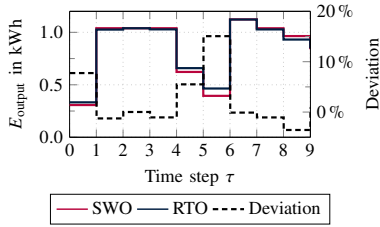
\begin{figure}[H]
    \centering
    \resizebox{!}{5cm}{\begin{tikzpicture}
\begin{axis} [ /pgf/number format/.cd,
        % use comma,
        1000 sep={.},
% legend style={at={(1,1)},anchor=north east}, 
% legend style={at={(0.5,-0.5)},anchor=north, legend columns=-1},
%width=0.49\linewidth, height=3cm, 
width=0.7\linewidth, height=5cm, 
   label style={font=\large},
xmin=0, xmax=9, 
xtick=data, 
% minor x tick num=1,  
ymin = 0, ymax = 1.2, 
minor y tick num=1,
xlabel={Time step $\tau$},
ylabel={$E_{\text{output}}$ in kWh},
yticklabel pos=left ,axis y line*=left, 
 ymajorgrids=true,
 xmajorgrids=true,
grid style = dotted, 
  y tick label style={
        /pgf/number format/.cd,
            fixed,
            fixed zerofill,
            precision=1,
        /tikz/.cd
    },
% tick style={draw=none}
]
\addplot[const plot, mark=, color=hsurot, very thick] table[y=Firststage, x = t] {figures/data_robustresults.txt};
\label{firststagelabel}
\addplot[mark=,mark size=1pt,const plot,color=hsublau, very thick]  table[y=RTO, x = t] {figures/data_robustresults.txt};
\label{rtolabel}
\end{axis}
\begin{axis} [/pgf/number format/.cd,
        % use comma,
        1000 sep={,},
legend style={at={(0.5,-0.35)},anchor=north, legend columns=-1},
width=0.7\linewidth, height=5cm, 
   label style={font=\large},
xmin=0, xmax=9, 
minor x tick num=0,  
ymin = -5, ymax = 20, 
yticklabel=\pgfmathprintnumber{\tick}\,$\%$,
visualization depends on={y \as \originalvalue},
minor y tick num=1,
xlabel={},
xtick=\empty,
ylabel={Deviation},
yticklabel pos=right ,axis y line*=right, 
 % ymajorgrids=true,
 % xmajorgrids=true,
% grid style = dotted, 
  y tick label style={
        /pgf/number format/.cd,
            fixed,
            fixed zerofill,
            precision=0,
        /tikz/.cd
    },
% tick style={draw=none}
]
\addlegendimage{/pgfplots/legend to name=firststagelabel, color  = hsurot, thick}
\addlegendentry{\ac{SWO}}
\addlegendimage{/pgfplots/legend to name=rtolabel, color  = hsublau, thick}
\addlegendentry{RTO}
\addplot[const plot, mark=, color=black, densely dashed, very thick]  table[y=DevPer, x = t] {figures/data_robustresults.txt};
\addlegendentry{Deviation}
\end{axis}
\end{tikzpicture}}
    \caption{Optimization results for entire horizon}
    \label{fig:results:robust}
\end{figure}

\autoref{fig:results:robust} shows the realized hydrogen output in each time step for the \ac{SWO} as well as the sum of all hydrogen outputs in the last \ac{RTO} time step of each \ac{SWO} time step $t_{\tau_j,9}$. Therein, it can be seen, that robust results have already been achieved through \ac{SWO}, with an average deviation of +2~\% between the optimized \textit{plan} and the realized \ac{SoS} across all time steps~$\tau$, encompassing a range of -3.5~\% to +15~\%. 
Although the \ac{SWO} generally yields favorable results, deviations exceeding +15~\% in periods of high \ac{RE} uncertainty underscore the importance of \ac{RTO}.

However, these inherent uncertainties in \ac{RE}, which resist further quantification and are dependend on the quality of the forecast, can be effectively addressed by \ac{RTO}, providing a valuable tool for appropriate responses.

The approach outlined in \autoref{sec:methodology} for implementing a two-stage optimization strategy is especially advantageous for hybrid energy resources that integrate both \ac{RE} sources and the electricity grid. Initially, the \ac{SWO} determines the required energy procurement from the energy market with manageable computational effort and a long-term planning horizon. The second optimization stage, \ac{RTO}, then utilizes the predetermined energy procurement from the energy market and the existing resource states established by the \ac{SWO} to respond to short-term fluctuations or deviations within optimized, defined framework conditions.

Utilizing one optimization model could not have accounted for continually updated forecasts as as this would have necessitated the continuous solving of the model for the entire horizon equal to $\tau$ at the availability of new forecasts. The computational times of one model at a high temporal resolution, e.g., 1.5 min, would have been too large to be useful (timely availability of the \textit{plan})~\cite{WKR+23}.
Such an approach is also unsuitable considering the high level of uncertainty in forecasts for the distant future, rendering a resolution of 1.5~min. unjustified, which emphasizes the application of the two-stage optimization approach presented.

In summary, the findings demonstrate the efficacy of the method presented in \autoref{sec:methodology} and affirm that the two-stage optimization approach can be successfully applied based on an existing static optimization model. This approach optimizes the procurement of energy from a spot market in parallel while facilitating the simultaneous integration of variable \ac{RE} sources. This is especially important given the increasing expansion of \ac{RE} and the imperative of decarbonization. %\cite{SCF+22}.

Moreover, it becomes evident that existing optimization models of \ac{SWO} can be adapted for this two-stage optimization approach without requiring significant adjustments (see \autoref{sec:method:challenges}). This facilitates the transition from pure planning to real time operation, accommodating short-term fluctuations.

% \subsection{Comparison of Method to Existing Optimization Approaches} \label{sec:vali:comparison}
% \todo[inline]{Unused potential benennen}
% Related work analyzed in \autoref{sec:relworks} 
\section{Discussion} \label{sec:discussion}
The method as well as the results of the evaluation, are discussed in this section.% regarding their interpretations, implications, and limitations.

%1. Interpretation
The presented research demonstrates the potential of transforming static optimization models into dynamic, \ac{RTO} models for managing flexible energy resources such as electrolyzers. Specifically, the work employs a two-stage optimization strategy that integrates both site-wide and real-time optimization techniques to improve operational flexibility. The method involves a hierarchical approach where long-term forecasts are used for \ac{SWO}, and short-term adjustments are made through \ac{RTO}. This allows the system to adapt dynamically to fluctuations in \ac{RE} generation and accurately reflect the actual behavior of the resource as well as the actual availability of \ac{RE}. Results from the case study of a system of electrolyzers validate the effectiveness of this approach. 

The proposed method ensures a stable power grid by matching energy consumption with \ac{RE} generation and demonstrates the practical benefits of integrating real-time adjustments into long-term planning. This dynamic capability is crucial for handling the variability and unpredictability associated with \ac{RE} sources.

%2. Implications
The implications of this study are relevant to the future of energy management and the integration of \ac{RE} sources. The two-stage optimization method developed in this work provides a robust framework for \ac{RTO} of flexible energy resources, such as electrolyzers. This method, by combining \ac{SWO} with \ac{RTO}, enhances the alignment of energy consumption with \ac{RE} generation. This alignment is important for sustainable energy practices, as it allows for more efficient use of \ac{RE} resources and minimizes waste and mismatch. This method makes it possible not to assume perfect foresight during the determination of a \ac{SoS} for subsequent control and rather continually use reliable values in \ac{RTO} in a short time gap to realization. This however leads to deviations from the initially generated \textit{plan} by \ac{SWO}. 

Such deviations, initially, are not problematic and are expected given the forecast uncertainties. In fact, if \ac{RTO} did not account for these deviations and instead strictly adhered to the \ac{SWO} \textit{plan} generated under the assumption of perfect foresight, unplanned adjustments would eventually become necessary. This is because either not all available \ac{RE} would be integrated or the grid supply would have to be adjusted in an unplanned and ad-hoc manner. Therefore, the ability of the \ac{RTO} to adapt to real-time data ensures more efficient and reliable integration of \ac{RE}, preventing the need for reactive measures and supporting overall grid stability.

%3. Limitations
Despite the promising outcomes, the work recognizes limitations. One notable limitation is the dependency on the accuracy of \ac{RE} generation forecasts. While the optimization models perform well with accurate forecasts, deviations from predicted values can influence the optimality of the \ac{SoS}. Additionally, the initial setup and computational requirements, although manageable, might present some challenges for large-scale implementation. The necessity for high-resolution data and computational resources could potentially impact the scalability of this approach. Future research should aim to enhance forecast accuracy, streamline computational processes, and explore decentralized optimization approaches to further improve the method's practicality and scalability.

Even though the case study showed the applicability of the method to a system of electrolyzers, the \ac{SWO} and \ac{RTO} of other types of flexible energy resources are also possible, as the underlying optimization model structure is generically applicable \cite{WRF23b, RWF23}.

\section{Conclusion}\label{sec:conclusion}
% \acresetall % Setze alle Akronymdefinitionen zurück
In this work, a novel two-stage optimization method was developed that bridges the gap between \ac{SWO} and \ac{RTO}, allowing for dynamic adjustments in response to the unpredictable nature of \ac{RE} sources. By leveraging existing static optimization models and transforming them for use in \ac{RTO}, a seamless integration of \ac{RE} into a system of electrolyzers is achieved, enhancing both grid reliability and operational efficiency. This method stands out for its adaptability, ensuring optimal resource utilization in near real-time and reducing the energy procurement costs from the intra-day market.

The case study on a system of electrolyzers, utilizing a hybrid energy supply from both the grid and a wind farm, validated the effectiveness of the method. The results demonstrated not only a more sustainable operation through the optimized use of renewable resources but also highlighted the potential for economic benefits as updated forecasts can be included in the operational planning. 

As the deployment of \ac{RE} sources continues to expand, the importance of such dynamic optimization methods will  increase. With higher shares of \ac{RE} in the power grid, the variability and unpredictability of the energy supply will become more pronounced. The proposed method is particularly relevant in this context, as it offers an adaptable solution to integrate increasing amounts of \ac{RE} while maintaining grid stability and optimizing resource utilization.

In conclusion, the developed two-stage optimization strategy offers a robust and flexible framework for energy management, paving the way for more sustainable and economically viable energy control of flexible energy resources in the face of growing \ac{RE} integration.

Future work can focus on the transfer of the method to decentralized optimization approaches employing  multi-agent systems, as they are inherently suited for the representation of resources under uncertainty and the handling thereof.
\appendix \label{appendix}
This appendix describes the mathematical modeling of \autoref{optmodel:sitewide}. For a detailed explanation please refer to \citet{WRF23b}. Additionally, the parameter set is presented.

\textit{Operational boundaries} are modeled as shown in \autoref{eq:power-boundaries-resource}~\cite{WRF23b}. 

\begin{align}
    P_{\text{min}, t} &\le  P_{t} \le  P_{\text{max}, t} \quad &\forall t \label{eq:power-boundaries-resource}
\end{align}

The piecewise linear approximation of the \textit{input-output relationship} is realized by means of binary variables $x_k$ for each segment $k \in \mathcal{K}$ and time step $\tau/t$ (\autoref{eq:pla-poweroutput}-\ref{eq_ub}). Only one segment can be active per time step (\autoref{eq:pla-sum1}). The total energy output  $D$ over the optimization horizon is set by \autoref{eq:target}. \cite{WRF23b}

\begin{align}
P_{\text{output},t} &= \sum \limits_{k\in \mathcal{K}} \left(a_{k} \cdot P_{\text{input}_k,t} + b_{k}\cdot x_{k,t} \right) \quad &\forall t \label{eq:pla-poweroutput}  \\
\text{lb}_{k} \cdot x_{k, t} &\le   P_{\text{input}_k, t}  &\forall t,k\label{eq_lb} \\
 P_{\text{input}_k, t} &\le  \text{ub}_{k}  \cdot x_{k, t} &\forall t,k\label{eq_ub} \\ 
 \quad \sum \limits_{k\in \mathcal{K}} x_{ k, t} &= 1 &\forall t\label{eq:pla-sum1} \\
  \Delta t \cdot  \sum \limits_{t\in\mathcal{T}} P_{t} &=  D \label{eq:target}
\end{align}

\textit{System states} $s \in \mathcal{S}$ are characterized by lower and upper flow limits (\autoref{eq:powerlimitsbystatemin} and~\ref{eq:powerlimitsbystatemax}), follower states $\mathcal{S}_{F,s}$  (\autoref{eq:statesequences}), holding durations of each state $s$ (\autoref{eq:minholdingduration} and~\ref{eq:maxholdingduration}), and ramp limits (\autoref{eq:rampmin} and~\ref{eq:rampmax}) \cite{WRF23b}.

\begin{align}
    \sum \limits_{s \in \mathcal{S}} x_{s, t} &= 1  &\forall t \label{eq:states-binary} \\
    P_{t} &\geq \sum \limits_{s \in \mathcal{S}} P_{\text{min}, s} \cdot x_{s, t} &\forall t>t_0 \label{eq:powerlimitsbystatemin} \\
    P_{t} &\leq \sum \limits_{s \in \mathcal{S}} P_{\text{max}, s} \cdot x_{s,t} &\forall t >t_0\label{eq:powerlimitsbystatemax} \\
    x_{t-1,s} - x_{t,s} &\le \sum \limits_{f \in \mathcal{S}_{F,s}} x_{f,t} &\forall s,t>t_0\label{eq:statesequences} \\
 t_{h, \text{min},s} \cdot  \left( x_{t,s} - x_{t-1,s}   \right) &\le \sum \limits_{\tau \in \mathcal{T}_h} x_{\tau,s}   &\quad \forall s, t>t_0 \label{eq:minholdingduration} \\
t_{h, \text{max},s} &\ge   \sum \limits_{\tau \in \mathcal{T}_h} x_{\tau,s}  &\forall s, t>t_0 \label{eq:maxholdingduration}   \\
\Delta t   \cdot \sum \limits_{s\in \mathcal{S}} \left(\text{ramp}_{\text{min},s} \cdot x_{s,t} \right) &\le \left|P_{t} - P_{t-1}\right| &\forall t>t_0 \label{eq:rampmin} \\
\Delta t   \cdot \sum \limits_{s\in \mathcal{S}} \left(\text{ramp}_{\text{max},s} \cdot x_{s,t} \right) &\ge  \left|P_{t} - P_{t-1}\right| &\forall t>t_0 \label{eq:rampmax}
\end{align}

The parameter set used for the case study has been derived from measurement data \cite{WaFa24}. \autoref{tab:parameters-pla} shows the parameters for the input-output relationship whereas \autoref{tab:parameters-states} shows parameters for system states and related constraints. These parameters are used for all electrolyzers.

\begin{table}[h]
\renewcommand{\arraystretch}{.7}
\centering
\caption{Parameters of the input-output relationship}
\label{tab:parameters-pla}
\begin{tabularx}{.85\linewidth}{p{2cm}*{4}{|Z}}
\toprule
Segment $k$  & 1 & 2 & 3 & 4\\
\cmidrule(lr){1-1} \cmidrule(lr){2-5}
lb, kW  & 0 & 0.6 & 1.2 & 1.8\\
ub, kW  & 0.6 & 1.2 & 1.8 & 2.4\\
$a_k$, kW/kW & 0.52 & 0.83 & 0.56 & 0.56\\
$b_k$, kW & -0.06 & -0.14 & 0.16 & 0.15\\
\bottomrule
\end{tabularx}
\end{table}

\begin{table}[h]
\renewcommand{\arraystretch}{.7}
\centering
\caption{System State Related Parameters}
\label{tab:parameters-states}
\begin{tabularx}{.75\linewidth}{p{2cm}*{3}{|Z}}
\toprule
State $s$ & 0 & 1 & 2 \\
% \midrule
\cmidrule(lr){1-1} \cmidrule(lr){2-4}
Name & off & stand-by & operation\\
$t_{h, \text{min},s}$ & 4 & 2 & 4\\
$t_{h, \text{max},s}$ & $\infty$ &  $\infty$ & $\infty$\\
$\mathcal{S}_{F,s}$ & \{2\} & \{0,2\}& \{0,1\}\\
$P_{\text{in., min},s}$, kW & 0 & 0.19 & 0.19\\
$P_{\text{in, max},s}$, kW& 0 & 0.19 & 2.4\\
$P_{\text{out, max},s}$, kW & 0 & 0 & 1.5\\
$\text{ramp}_{\text{min},s}$, kW/h & 0 & 0 & 0 \\
$\text{ramp}_{\text{max},s}$, kW/h & 25000 & 3456 &3456\\
\bottomrule
\end{tabularx}
\end{table}

\bibliographystyle{IEEEtranN}
%\footnotesize{\bibliography{literature}}
\footnotesize{% Generated by IEEEtranN.bst, version: 1.14 (2015/08/26)

}
\end{document}